\documentclass[10pt,a4paper]{article}
\usepackage{amssymb}
\usepackage{amscd}
\usepackage{latexsym}
\usepackage{amsmath}
\usepackage{amsfonts}
\usepackage{amscd}

\begin{document}

\newtheorem{theorem}{Theorem}[section]
\newtheorem{remark}[theorem]{Remark}
\newtheorem{mtheorem}[theorem]{Main Theorem}
\newtheorem{bbtheo}[theorem]{The Strong Black Box}
\newtheorem{observation}[theorem]{Observation}
\newtheorem{proposition}[theorem]{Proposition}
\newtheorem{lemma}[theorem]{Lemma}
\newtheorem{testlemma}[theorem]{Test Lemma}
\newtheorem{mlemma}[theorem]{Main Lemma}
\newtheorem{note}[theorem]{{\bf Note}}
\newtheorem{steplemma}[theorem]{Step Lemma}
\newtheorem{corollary}[theorem]{Corollary}
\newtheorem{notation}[theorem]{Notation}
\newtheorem{example}[theorem]{Example}
\newtheorem{definition}[theorem]{Definition}

\renewcommand{\labelenumi}{(\roman{enumi})}
\def\Pf{\smallskip\goodbreak{\sl Proof. }}

\def\Fin{\mathop{\rm Fin}\nolimits}
\def\br{\mathop{\rm br}\nolimits}
\def\fin{\mathop{\rm fin}\nolimits}
\def\Ann{\mathop{\rm Ann}\nolimits}
\def\Aut{\mathop{\rm Aut}\nolimits}
\def\End{\mathop{\rm End}\nolimits}
\def\bfb{\mathop{\rm\bf b}\nolimits}
\def\bfi{\mathop{\rm\bf i}\nolimits}
\def\bfj{\mathop{\rm\bf j}\nolimits}
\def\df{{\rm df}}
\def\bfk{\mathop{\rm\bf k}\nolimits}
\def\bEnd{\mathop{\rm\bf End}\nolimits}
\def\iso{\mathop{\rm Iso}\nolimits}
\def\id{\mathop{\rm id}\nolimits}
\def\Ext{\mathop{\rm Ext}\nolimits}
\def\Ines{\mathop{\rm Ines}\nolimits}
\def\Hom{\mathop{\rm Hom}\nolimits}
\def\bHom{\mathop{\rm\bf Hom}\nolimits}
\def\Rk{ R_\k-\mathop{\bf Mod}}
\def\Rn{ R_n-\mathop{\bf Mod}}
\def\map{\mathop{\rm map}\nolimits}
\def\cf{\mathop{\rm cf}\nolimits}
\def\top{\mathop{\rm top}\nolimits}
\def\Ker{\mathop{\rm Ker}\nolimits}
\def\Bext{\mathop{\rm Bext}\nolimits}
\def\Br{\mathop{\rm Br}\nolimits}
\def\dom{\mathop{\rm Dom}\nolimits}
\def\min{\mathop{\rm min}\nolimits}
\def\im{\mathop{\rm Im}\nolimits}
\def\max{\mathop{\rm max}\nolimits}
\def\rk{\mathop{\rm rk}}
\def\Diam{\diamondsuit}
\def\Z{{\mathbb Z}}
\def\Q{{\mathbb Q}}
\def\N{{\mathbb N}}
\def\bQ{{\bf Q}}
\def\bF{{\bf F}}
\def\bX{{\bf X}}
\def\bY{{\bf Y}}
\def\bHom{{\bf Hom}}
\def\bEnd{{\bf End}}
\def\bS{{\mathbb S}}
\def\AA{{\cal A}}
\def\BB{{\cal B}}
\def\CC{{\cal C}}
\def\DD{{\cal D}}
\def\TT{{\cal T}}
\def\FF{{\cal F}}
\def\GG{{\cal G}}
\def\PP{{\cal P}}
\def\SS{{\cal S}}
\def\XX{{\cal X}}
\def\YY{{\cal Y}}
\def\fS{{\mathfrak S}}
\def\fH{{\mathfrak H}}
\def\fU{{\mathfrak U}}
\def\fW{{\mathfrak W}}
\def\fK{{\mathfrak K}}
\def\PT{{\mathfrak{PT}}}
\def\T{{\mathfrak{T}}}
\def\fX{{\mathfrak X}}
\def\fP{{\mathfrak P}}
\def\X{{\mathfrak X}}
\def\Y{{\mathfrak Y}}
\def\F{{\mathfrak F}}
\def\C{{\mathfrak C}}
\def\B{{\mathfrak B}}
\def\J{{\mathfrak J}}
\def\fN{{\mathfrak N}}
\def\fM{{\mathfrak M}}
\def\Fk{{\F_\k}}
\def\bar{\overline }
\def\Bbar{\bar B}
\def\Cbar{\bar C}
\def\Pbar{\bar P}
\def\etabar{\bar \eta}
\def\Tbar{\bar T}
\def\fbar{\bar f}
\def\nubar{\bar \nu}
\def\rhobar{\bar \rho}
\def\Abar{\bar A}
\def\a{\alpha}
\def\b{\beta}
\def\g{\gamma}
\def\w{\omega}
\def\e{\varepsilon}
\def\o{\omega}
\def\va{\varphi}
\def\k{\kappa}
\def\m{\mu}
\def\n{\nu}
\def\r{\rho}
\def\f{\phi}
\def\hv{\widehat\v}
\def\hF{\widehat F}
\def\v{\varphi}
\def\s{\sigma}
\def\l{\lambda}
\def\lo{\lambda^{\aln}}
\def\d{\delta}
\def\z{\zeta}
\def\th{\theta}
\def\t{\tau}
\def\ale{\aleph_1}
\def\aln{\aleph_0}
\def\Cont{2^{\aln}}
\def\nld{{}^{ n \downarrow }\l}
\def\n+1d{{}^{ n+1 \downarrow }\l}
\def\hsupp#1{[[\,#1\,]]}
\def\size#1{\left|\,#1\,\right|}
\def\Binfhat{\widehat {B_{\infty}}}
\def\Zhat{\widehat \Z}
\def\Mhat{\widehat M}
\def\Rhat{\widehat R}
\def\Phat{\widehat P}
\def\Fhat{\widehat F}
\def\fhat{\widehat f}
\def\Ahat{\widehat A}
\def\Chat{\widehat C}
\def\Ghat{\widehat G}
\def\Bhat{\widehat B}
\def\Btilde{\widetilde B}
\def\Ftilde{\widetilde F}
\def\restr{\mathop{\upharpoonright}}
\def\to{\rightarrow}
\def\arr{\longrightarrow}
\def\LA{\langle}
\def\RA{\rangle}
\newcommand{\norm}[1]{\text{$\parallel\! #1 \!\parallel$}}
\newcommand{\supp}[1]{\text{$\left[ \, #1\, \right]$}}
\def\set#1{\left\{\,#1\,\right\}}
\newcommand{\mb}{\mathbf}
\newcommand{\wt}{\widetilde}
\newcommand{\card}[1]{\mbox{$\left| #1 \right|$}}
\newcommand{\union}{\bigcup}
\newcommand{\inters}{\bigcap}
\newcommand{\ER}{{\rm E}}
\def\Proof{{\sl Proof.}\quad}
\def\fine{\ \black\vskip.4truecm}
\def\black{\ {\hbox{\vrule width 4pt height 4pt depth
0pt}}}
\def\fine{\ \black\vskip.4truecm}
\long\def\alert#1{\smallskip\line{\hskip\parindent\vrule%
\vbox{\advance\hsize-2\parindent\hrule\smallskip\parindent.4\parindent%
\narrower\noindent#1\smallskip\hrule}\vrule\hfill}\smallskip}

\title{Completely Inert Subgroups of Abelian Groups}
\author{Andrey R. Chekhlov \\Faculty of Mathematics and Mechanics, Section of Algebra, \\Tomsk State University, Tomsk 634050, Russia\\{\small e-mails: cheklov@math.tsu.ru, a.r.che@yandex.ru}
\\and\\
Peter V. Danchev \\Institute of Mathematics and Informatics, Section of Algebra, \\Bulgarian Academy of Sciences, Sofia 1113, Bulgaria\\{\small e-mails: danchev@math.bas.bg, pvdanchev@yahoo.com}}
\maketitle

\centerline{(To Patrick W. Keef on the occasion of his {\bf 70}th birthday)}

\medskip
\medskip

\begin{abstract}{We define and study in-depth the so-called {\it completely inert} and {\it uniformly completely inert} subgroups of Abelian groups. We curiously show that a subgroup is completely inert exactly when it is characteristically inert. Moreover, we prove that a subgroup is uniformly completely inert precisely when it is uniformly characteristically inert. These two statements somewhat strengthen recent results due to Goldsmith-Salce established for totally inert subgroups in J. Commut. Algebra (2025).

Some other closely relevant things are obtained as well.}
\end{abstract}

\section{Some Fundamentals}

Throughout the current brief paper, all our groups are {\it additively} written and {\it Abelian}. Our notation and terminology are mainly standard and follow those from \cite{Kap}. In fact, recall the standard concepts that an arbitrary subgroup $F$ of a group $G$ is said to be {\it fully invariant} provided $\phi(F)\subseteq F$ for any endomorphism $\phi$ of $G$, and an arbitrary subgroup $C$ of $G$ is said to be a {\it characteristic} subgroup provided $\psi(F)\subseteq F$ for any automorphism $\psi$ of $G$. Moreover, in \cite{C}, an arbitrary subgroup $S$ of $G$ is said to be a {\it strongly invariant} subgroup provided $f(S)\subseteq S$ for any homomorphism $f: S\to G$.

\medskip

It is obvious that strongly invariant subgroups are fully invariant subgroups, and the later are always characteristic, while both the reverse fail in general.

\medskip

Generalizing these two notions, it is well known that a subgroup $N$ of a group $G$ is called {\it fully inert} provided $(\phi(N)+N)/N$ is finite for all endomorphisms $\phi$ of $G$, and is called {\it characteristically inert} provided $(\psi(N)+N)/N$ is finite for all automorphisms $\psi$ of $G$ (for the latter see \cite{CDG2}). Besides, in \cite{BC}, a subgroup $N$ of $G$ is called {\it strongly inert} provided $(f(N)+N)/N$ is finite for all homomorphisms $f: N\to G$. If, in both cases, the cardinalities of the finite quotients $(\phi(N)+N)/N$ and $(\psi(N)+N)/N$ are bounded by some fixed positive integer, the subgroup $N$ is termed {\it uniformly fully inert} and {\it uniformly characteristically inert}, respectively. The description of these subgroups as being commensurable with fully invariant and, respectively, with characteristic subgroups can be found in \cite{CD}.

\medskip

Apparently, strong inertness yields full inertness yields characteristic inertness, whereas the opposite are both generally untrue.

\medskip

Further, in \cite{GS}, the new concept of a totally inert subgroup of a group was introduced as follows: A subgroup $T$ of an arbitrary group $G$ is called {\it totally inert} provided the intersection $T\cap \phi(T)$ has finite index in both $T$ and $\phi(T)$ for any non-zero endomorphism $\phi$ of $G$.

\medskip

Clearly, total inertness implies full inertness.

\medskip

Imitating the "uniformly" property presented above, it is reasonably natural to ask what is the behavior of {\it uniformly totally inert} subgroups defined analogously as follows: The subgroup $T$ has an intersection $T\cap \phi(T)$ bounded by a fixed positive integer in $T$ and $\phi(T)$ for any non-zero endomorphism $\phi$ of $G$.

\medskip

However, this does {\it not} give nothing new as the next arguments illustrate -- indeed, we claim that these are only the rational torsion-free group $\mathbb{Q}$ and the quasi-cyclic $p$-group $\mathbb{Z}(p^{\infty})$. In fact, looking for infinite subgroups $H$ of a group $G$ different from $\mathbb{Q}$ and $\mathbb{Z}(p^{\infty})$, respectively, which are uniformly totally inert, from \cite[Corollary 2.3]{GS} it follows that, in order to admit infinite totally inert subgroups, $G$ must be torsion-free reduced and indecomposable. Furthermore, for such a group $G$, a subgroup $H \neq \{0\}$ is infinite, but and it cannot be uniformly totally inert, because there is a prime $p$ such that $H/pH \neq \{0\}$ (noticing that $H$ is not divisible), so $H/pH$ has cardinality at least $p$ and, therefore, $H/p^nH$ has cardinality at least $p^n$, whence $H$ cannot be uniformly totally inert, as expected.

\medskip

In order to strengthen this, we come to the following basic tool (see the initial version in \cite{CDK} as well).

\begin{definition} A subgroup $C$ of an arbitrary group $G$ is called {\it completely inert} provided the intersection $C\cap \psi(C)$ has finite index in both $C$ and $\psi(C)$ for any automorphism $\psi$ of $G$.
\end{definition}

Evidently, total inertness forces complete inertness, but the reciprocal implication is manifestly non-true as we will illustrate in the sequel.

\medskip

As above, a reasonably logical question is to ask what happens with {\it uniformly completely inert subgroups} defined by analogy thus: The subgroup $C$ has an intersection $C\cap \psi(C)$ bounded by a fixed positive integer in $C$ and $\psi(C)$ for any automorphism $\psi$ of $G$.

\medskip

And so, the objective of this article is to give a systematic exploration of the so-introduced concept of complete inertness by finding its crucial properties and comparing them with these of the defined above total inertness.

\section{Main Results}

Before establishing the principal achievements, our preliminaries here are the following. The first technicality can easily be established by analogy with \cite[Proposition 2.1]{GS}.

\begin{lemma}\label{1} A subgroup commensurable with a completely inert subgroup is completely inert.
\end{lemma}

\Pf It is straightforwardly analogous to Proposition 2.1 from \cite{GS}.
\fine

We now arrive at the following quite surprising assertion.

\begin{proposition}\label{2} A subgroup $C$ of a group $G$ is completely inert subgroup if, and only if, it is characteristically inert.
\end{proposition}

\Pf Necessity is evident, so we omit the arguments. As for sufficiency, by assumption, the quotient $\psi^{-1}(C)/(C\cap\psi^{-1}(C))$ is finite for each $\psi\in\mathrm{Aut}(G)$. We claim that the factor-group $C/(C\cap\psi(C))$ is likewise finite as being an isomorphic image of $\psi^{-1}(C)/(C\cap\psi^{-1}(C))$ under the action of $\psi$.

In fact, mapping $$\overline{\psi}: \psi^{-1}(c)+(C\cap\psi^{-1}(C))\mapsto c+(\psi(C)\cap C),$$ one concludes that it is an isomorphism between $\psi^{-1}(C)/(C\cap\psi^{-1}(C))$ and $C/(C\cap\psi(C))$.
Indeed, $\psi$ maps $C\cap\psi^{-1}(C)$ into $\psi(C)\cap C$, so that $\overline{\psi}$ is obviously a well-defined homomorphism. Furthermore, if $c=\psi(c_1)$ for some $c,c_1\in C$, then $c_1=\psi^{-1}(c)\in C$,
and hence $\overline{\psi}$ is an injection. It is also routinely seen that $\overline{\psi}$ is a surjection. Whence, $\overline{\psi}$ is an isomorphism, as asserted.
\fine

Our next pivotal instrument is the following.

\begin{proposition}\label{3} A subgroup $U$ of a group $G$ is uniformly completely inert if, and only if, it is uniformly characteristically inert.
\end{proposition}

\Pf The same idea is workable since, by assumption, there is $n\in \mathbb{N}$ such that $|\psi^{-1}(C)/(C\cap\psi^{-1}(C))|\leq n$ for each $\psi\in\mathrm{Aut}(G)$, and so we claim that $|C/(C\cap\psi(C))|\leq n$
as being an isomorphic image of $\psi^{-1}(C)/(C\cap\psi^{-1}(C))$ acting by $\psi$, as demonstrated above.
\fine

Recall that a group $G$ is said to have {\it unit sum number} $\mathrm{usn}(G)=n\in\mathbb{N}$ if each endomorphism of $G$ is a sum of $\leq n$ automorphisms of $G$.

\medskip

As the following lemma illustrates, if $\mathrm{usn}(G)=n$, then every uniformly characteristically inert subgroup of $G$ is uniformly fully inert in $G$, hopefully bounded by another fixed positive integer.

\begin{lemma}\label{0.2}
(1) If $\mathrm{usn}(G)=n$, then every uniformly characteristically inert subgroup $C$ of $G$ (bounded by a fixed integer $k>0$) is uniformly fully inert in $G$ (bounded by an integer $\leq nk$).

(2) If $C$ is uniformly characteristically inert in $G$ and $H\sim C$, then $H$ is uniformly characteristically inert in $G$.
\end{lemma}

\Pf (1) Letting $\varphi$ be an arbitrary endomorphism of $G$, then there are $\alpha_1,\dots,\alpha_n\in\mathrm{Aut}\,G$ with $\varphi=\alpha_1+\dots+\alpha_n$. But then we obtain
$$(C+\varphi(C))/C\leq (\alpha_1(C)+C)/C+\dots+(\alpha_n(C)+C)/C,$$ and so $$|(C+\varphi(C))/C|\leq |(\alpha_1(C)+C)/C|+\dots+|(\alpha_n(C)+C)/C|\leq nk,$$ as claimed.

(2) It can be verified in two different ways:

(a) It must be shown that the order of the factor-group $\alpha(H)/(\alpha(H)\cap H)$ is finite and the same for all $\alpha\in\mathrm{Aut}\,G$. To that end, since $|H/(H\cap C)|,|C/(H\cap C)|\leq n$ for some
$n$, one inspects that $$|\alpha(H)/(\alpha(H)\cap\alpha(C))|,|\alpha(C)/(\alpha(H)\cap\alpha(C))|\leq n$$ for all $\alpha\in\mathrm{Aut}\,G$. So, $\alpha(C)\sim\alpha(H)$, but as $\alpha(C)\sim C$, in view of $C\sim H$
we can get $\alpha(H)\sim H$ for any $\alpha\in\mathrm{Aut}\,G$.

(b) Knowing that each uniformly characteristically inert subgroup is commensurable with some characteristic subgroup (see \cite[Corollary 1.9]{DDR}), we may deduce that $H$ is commensurable with some characteristic subgroup,
and consequently is uniformly characteristically inert, as asserted.
\fine

We now proceed by proving a series of technicalities.

\begin{lemma} Let $H$ be an infinite uniformly characteristically inert subgroup with bounded index $k\geq 1$ of the group $G=\bigoplus_{i\in I} G_i$, where $G_i\cong G_j$ for all $i,j\in I$ and the index set
$I$ is infinite. If $\pi_i$ denotes the canonical projection of $G$ onto $G_i$, then $H\leq \bigoplus_{i\in I}\pi_i(H)$ with bounded index $3k$.
\end{lemma}

\Pf Recall that $\mathrm{usn}(G)=3$ (see, e.g., \cite{CDG2}). Now, Lemma~\ref{0.2} tells us that $H$ is, actually, uniformly fully inert. If, however, we assume the contrary that $\bigl|\bigl(\bigoplus_{i\in I}\pi_i(H)\bigr)/H\bigr|>3k$, then there will exist $i_1,\dots,i_n\in I$ such that, for the sum $\psi=\pi_{i_1}+\dots+\pi_{i_n}$, we have $$|(\psi(H)+H)/H|=|(\pi_{i_1}(H)+\dots+\pi_{i_n}(H))/H|>3k,$$ contradicting to Lemma~\ref{0.2}, as suspected.
\fine

Our next technical statement is the following known fact.

\begin{lemma}~\label{0.4} (\cite[Lemma 2.2]{CDG1}) Let $H$ be a fully inert subgroup of the group $G=\bigoplus_{i\in I} G_i$, where the index set $I$ is infinite, and let each $\pi_i$ denote the canonical projection from $G$ onto
$G_i$. Then, $H$ is commensurable with $\bigoplus_{i\in I}\pi_i(H)$, the images $\pi_i(H)$ are fully inert in $G_i$, and almost all $\pi_i(H)$ are fully invariant in $G_i$. Furthermore, there is a finite subset
$S\subset I$, such that $\bigoplus_{i\in I\setminus S}\pi_i(H)$ is fully invariant in $\bigoplus_{i\in I\setminus S}G_i$.
\end{lemma}

We are, thereby, ready to attack the following.

\begin{proposition}\label{27} Let $H$ be a fully inert subgroup of the group $G=\bigoplus_{i\in I} G_i$, where $G_i\cong G_j$ for every $i,j\in I$ and the index set $I$ is infinite. Then, $H$ is commensurable with some fully invariant subgroup of $G$.
\end{proposition}

\Pf Invoking Lemma~\ref{0.4}, with no loss of generality we can assume that $H=\bigoplus_{i\in I}\pi_i(H)$. For $j\in I\setminus S$, denoting $F_j=\pi_j(H)$, we derive that $\bigoplus_{j\in I\setminus S}F_j$
is fully invariant in $\bigoplus_{j\in I\setminus S}G_j$, where it is clear that all $F_j$ are isomorphic ($j\in I\setminus S$).

After that, in each $G_i$, where $i\in S$, there exists a fully invariant subgroup $F_i\cong F_j$ with $j\in I\setminus S$. Thus, $F=\bigoplus_{i\in I}F_i$ is a fully invariant subgroup of $G$. It next suffices to show that, for each $i\in S$, the subgroup $F_i$ is commensurable with $\pi_i(H)$. To this aim, let $i$ be some index in $S$ and $\varphi\in\mathrm{E}(G)$ such that $\varphi(G_i)=G_{j_0}$, where $j_0$ is a fixed index from $I\setminus S$ with $\varphi(G_k)=\{0\}$ whenever $k\neq i$.
So, one checks that $$H+\varphi(H)=\bigl(\bigoplus_{s\in S}\pi_s(H)\bigr)\oplus(F_{j_0}+F'_{j_0})\oplus\bigl(\bigoplus_{l\in I\setminus (S\cup\{j_0\})}F_l\bigr),$$
where $F_{j_0}'=\varphi\pi_i(H)$, whence $F_{j_0}+F'_{j_0}\sim F_{j_0}$.

Let us now $\psi\in\mathrm{E}(G)$ such that $\psi(G_{j_0})=G_i$ and $\varphi(G_k)=\{0\}$ whenever $k\neq j_0$. Therefore, one infers that
$$H+\psi(H)=\bigl(\bigoplus_{s\in S\setminus \{i\}}\pi_s(H)\bigr)\oplus(F'_i+\pi_i(H))\oplus\bigl(\bigoplus_{l\in I\setminus S}F_l\bigr),$$
where $F_i'=\psi(F_{j_0})=F_i$, whence $F_i+\pi_i(H)\sim\pi_i(H)$.

But since $F_{j_0}'\cong\pi_i(H)$, $F_{j_0}\cong F_i$ and these isomorphisms both induce endomorphisms of the group $G$, we then arrive at the relations $F_i+\pi_i(H)\sim F_i$ and $F_i+\pi_i(H)\sim\pi_i(H)$, i.e., $\pi_i(H)\sim F_i$, as asked for.
\fine

\begin{lemma} If $G$ is a group in which each fully inert subgroup is uniformly fully inert and $A$ is a direct summand in $G$, then in $A$ each fully inert subgroup is uniformly fully inert.
\end{lemma}

\Pf Write $G=A\oplus B$, $H\leq A$ is fully inert in $A$ and set $F:=\mathrm{Hom}(A,B)H$. Since $F$ is fully invariant in $B$, then it is not too hard to establish that the subgroup $H\oplus F$ is fully inert in $G$, and so by assumption $H\oplus F$ is uniformly fully inert in $G$. However, every endomorphism of $A$ can be considered as an endomorphism of $G$, and therefore $H$ is uniformly fully inert in $A$, as pursued.
\fine

The next technical assertion is well-known.

\begin{lemma}~\label{0.7} (\cite[Lemma 7]{C17}) Let $G=\bigoplus_{i\in I} G_i$, and let $\pi_i: G\to G_i$ be the corresponding projections. If $H$ is a fully inert subgroup of $G$, then $H$ is of finite index in the subgroup
$\bigoplus_{i\in I}\pi_i(H)$.
\end{lemma}

It is worthy of noticing that, if $G:=A^n$ for some natural number $n$ and a group $A$, then not every fully inert subgroup of $G$ is commensurable with a fully invariant subgroup (see, for example, \cite{C22}).

\medskip

We finish the series of technical claims with the last one.

\begin{lemma}\label{210} Let $G=A_1\oplus\dots\oplus A_n$, $A_i\cong A$ for $i=1,\dots,n$, where $A$ is a non-zero group, and suppose in $A$ each fully inert subgroup is commensurable with a fully invariant subgroup. Then, each fully inert subgroup $H$ of $G$ is commensurable with a fully invariant subgroup of $G$. In particular, the subgroup $H$ is uniformly fully inert.
\end{lemma}

\Pf Utilizing Lemma~\ref{0.7}, it can be assumed that $H=\pi_1(H)\oplus\dots\oplus\pi_n(H)$. It is easy to see that each $\pi_i(H)$ is fully inert in $A_i$ for $i=1,\dots,n$. So, each $\pi_i(H)$ is commensurable with a
fully invariant subgroup $F_i$ of $A_i$. Likewise, it is evident that each $F_i$ is commensurable with $\phi(F_j)$ for each isomorphism $\phi: A_j\to A_i$. Thus, it can be supposed that $F_i\cong F_j$ for all $i,j=1,\dots,n$.
Consequently, in this case, $F=F_1\oplus\dots\oplus F_n$ is a fully invariant subgroup of $G$ which is commensurable with $H$, as required.
\fine

The following construction is worthwhile, although it was documented in \cite[Example 4.1]{CDG2} that there is a group with a characteristically inert subgroup which is {\it not} fully inert. Thus, Proposition~\ref{2} yields that there is a group with completely inert subgroup which is {\it not} fully inert and so, manifestly, {\it not} totally inert.

\begin{example} There exists a group with a completely inert subgroup which is {\it not} totally inert.
\end{example}

\Pf Let $G$ be a torsion-free group of infinite rank such that $\mathrm{E}(G)\cong\mathbb{Z}$, the group of integers. Therefore, each its subgroup is characteristic and hence is completely inert. But $H\leq G$ will be totally inert exactly when the factor-group $H/nH$ is finite for all non-zero integers $n$, that is, if it is a narrow subgroup in terms of \cite{GS}. But, because $G$ must have infinite rank, then by consulting with \cite[Theorem 2.8]{GS} in $G$
there will exist subgroups that are not totally inert, as promised.
\fine

\begin{remark} Concerning the preceding example, note that if $G$ is a torsion-free group each subgroup of which is completely inert and $pG=G$ for at least one prime $p$, then the rank of $G$ is necessarily finite. Indeed, if the rank of $G$ is infinite, then $G$ will have a free subgroup $H$ of infinite rank with the property $H\nsim pH$.

Likewise, there is an abundance of completely inert subgroups which are not totally inert: in fact, such are all infinite characteristic subgroups of decomposable groups $A$ (e.g., the subgroups $nA$ and $(p^mA)[n]$ for some naturals $m,n$) as well as, if a group $G$ is torsion-free, then such subgroups are $G(t)$ and $G(\chi)$, where $t$ is the type and $\chi$ is the characteristic of $G$.
\end{remark}

In the other vein, since each totally inert subgroup of a decomposable group is finite (see \cite[Proposition 2.2]{GS}), it must be that each infinite characteristic subgroup of a decomposable group will be completely inert but {\it not} totally inert subgroup.

\medskip

We now proceed by establishing a few more statements in this directory.

\begin{lemma}\label{0001} Suppose in a group $G$ every subgroup is completely inert. The following two points hold:

(1) if $G=A\oplus B$, then $f(B)$ is finite for every homomorphism $f: B\to A$;

(2) if $G$ is torsion-free, then $G$ is indecomposable.
\end{lemma}

\Pf (1) It is evident since $\varphi=1_G+f\in\mathrm{Aut}\,G$ and $\varphi(B)=f(B)\oplus B\sim B$.

(2) Let $H=\langle a+b\rangle$, where $0\neq a\in A$, $0\neq b\in B$. If $f=1_A\oplus(-1_B)$, then $f(H)=\langle a-b\rangle$ and $H\cap f(H)=\{0\}$, and so $H\nsim f(H)$, as required.
\fine

Thus, from this lemma, it immediately follows that in any divisible group $D$ each subgroup is completely inert if, and only if, either $D\cong\mathbb{Q}$, the group of all rational numbers, or $D$ is torsion and the rank of each its $p$-component is $\leq 1$.

\begin{lemma}\label{0002} The following three items hold:

(1) If $G$ is a reduced $p$-group in which each subgroup is completely inert, then $G$ is finite.

(2) In a non-reduced $p$-group $G$ all subgroups are completely inert if, and only, if $G$ is a direct sum
of a quasi-cyclic group and a finite group.

(3) In a torsion group all subgroups are completely inert if, and only if, almost all its $p$-components are co-cyclic and the non co-cyclic $p$-components are either finite or a direct sum of a quasi-cyclic group and a finite group.
\end{lemma}

\Pf (1) From Lemma~\ref{0001}, it follows that the Ulm-Kaplansky invariant $f_n(G)$ is finite for all $n<\omega$. So, if we assume the contrary that $G$ is infinite, and thus unbounded, then its basic subgroup, say $B$, also is unbounded. Choose in $B$ such a direct summand $H=\bigoplus_{i=1}^{\infty}\langle a_i\rangle$ that in the complementary direct summand there exists a direct summand $F=\bigoplus_{i=1}^{\infty}\langle b_i\rangle$ with ${\rm order}(b_i)=p^{n_i}$ and all $n_i\geq 1$.
Therefore,
$\overline{H}=(H+pG)/pG$ is a direct summand in
$G/pG$, so that there exists a homomorphism $f: G/pG\to G[p]$ such that $f((a_i+pG)/pG)=p^{n_i-1}b_i$ and $f$ acts as the zero homomorphism on the additional direct summand of $\overline{H}$.

Now, if $\pi:G\to G/pG$ is the usual canonical surjection and $\varphi=f\circ\pi$, then one checks that $\varphi^2=0$ and thus $\psi=1+\varphi\in\mathrm{Aut}\,G$. But, as $H\cap\psi(H)=pH$, we infer $H\cap\psi(H)\nsim H$, proving that $G$ is finite, as promised.

Furthermore, items (2) and (3) follows directly from Lemma~\ref{0001} and item (1).
\fine

\begin{lemma}\label{0003} If $G=\mathbb{Q}\oplus B$ is a group, where $\mathbb{Q}$ is the rational group, then in $G$ each subgroup is completely inert if, and only, if $B$ is finite.
\end{lemma}

\Pf {\bf Necessity.} An appeal to Lemma~\ref{0001} riches us that $B$ is torsion, almost all $p$-components of $B$
are cyclic and the non-cyclic $p$-components are finite. If, in a way of contradiction, we assume for a moment that $B$ is infinite (i.e., an infinite number of its $p$-components are non-zero), then in $\mathbb{Q}$ there will exist a subgroup, $X$ say, with an infinite image of the existing homomorphism $X\to B$, that this can not be happen, contrary to our assumption.

{\bf Sufficiency.} It is obvious, so the arguments are removed voluntarily.
\fine

Given a reduced torsion-free group $G$, we shall denote by $R(G)$ the maximal subring of the field of rational numbers $\mathbb{Q}$ contained in $\mathrm{E}(G)$, which is generated by the $1$ and the inverses of the prime numbers $p$ for which $G = pG$.

\begin{lemma}\label{0004} Let $G$ be a torsion-free group such all endomorphism are injections. Then, in $G$ each subgroup is completely inert if, and only if, $\mathrm{Aut}\,G=U(R(G))$ and, moreover, if $R(G)\neq\mathbb{Z}$, then the rank of $G$ is finite.
\end{lemma}

\Pf {\bf Necessity.} Suppose that $H$ is a pure subgroup of rank $1$, and $\varphi\in\mathrm{Aut}\,G$. Since
$H\sim\varphi(H)$, we derive $H=\varphi(H)$, so $\varphi\!\upharpoonright\! H\in R(H)$ and hence $\varphi\!\upharpoonright\! H=m/n$, where $m$ and $n$ are mutually simple integers. However, as $(n\varphi-m)H=\{0\}$, we deduce $n\varphi=m$. Therefore, $nG=G$ whence $\varphi=m/n\in R(G)$.

{\bf Sufficiency.} If $R(G)=\mathbb{Z}$, then one sees that each subgroup is totally inert. If $pG=G$ for some prime
$p$, then under assumption $G$ possesses finite rank, and so each $H\leq G$ will too have finite rank, whence
$H\sim\varphi(H)$ for every $\varphi\in R(G)$.
\fine

Let us now recollect the following technicality necessary for our further presentation.

\begin{lemma}\label{fuch1}~(\cite[\S 8, exercise 5]{F}) If $H$ is a subgroup of $G=B\oplus C$, then $H$ is the sub-direct sum of the groups $B\cap(H+C)$ and $(B+H)\cap C$.
\end{lemma}

We now intend to prove the following statement.

\begin{proposition}\label{0006} In the splitting group $G=T\oplus R$, where $T=T(G)$, each subgroup is completely inert if, and only if, in both $T$ and $R$ each subgroup is completely inert and, moreover, the image of every homomorphism $f: S\to T$ is finite for every subgroup $S\leq R$.
\end{proposition}

\Pf {\bf Necessity.} Assume that $f(S)$ is infinite for some $S\leq R$. Letting $H:=\{f(x)+x\,|\,x\in S\}$, then for
$\varphi:=1_T\oplus (-1_R)$ we have
$$T'=2f(S)=\{2f(x)\,|\,x\in S\}\leq H+\varphi(H)\sim H.$$ So, $T'$ is finite, because $H\cap T=\{0\}$. However, since $2T_q(f(S))=T_q(f(S))$ for any prime $q>2$, $T_2(f(S))$ has to be infinite, and so it is non-reduced. But, in this case, $2(f(S))$ must be infinite, and this contradiction guarantees that
$f(S)$ is finite.

{\bf Sufficiency.} Assume $H\leq G$. Then, Lemma~\ref{fuch1} informs us that $H$ is a sub-direct sum of the groups
$T_1$ and $R_1$, where $T_1=T\cap(H+R)$ and $R_1=(T+H)\cap R$. Consequently, if $T_0:=T_1\cap H$ and $R_0:=H\cap R_1$, then we know with help of the property of sub-direct sums that $T_1/T_0\cong R_1/R_0$.

If, for a moment, the factor-group $T_1/T_0$ is non-reduced, then given the structure of the group $T$, we obtain that non-reduced must be the group $T$ as well; thus, there is a homomorphism $R_1\to T$ with infinite image, that is wrong.

If, however, $T_1/T_0$ is reduced, then each of its $p$-components is a factor-group of either a cyclic group or of a finite group, and thus there exists a homomorphism $T_1/T_0\to T_1$ with non-zero image of each non-zero $p$-component of $T_1/T_0$. Hence, if $T_1/T_0$ is infinite, then there will exist a homomorphism $R_1\to T$ with infinite image, that is untrue.

Therefore, $T_1/T_0$ has to be finite, and so finite is $R_1/R_0$ too. Furthermore, since $\theta(H)=R_1\sim R_0$,
where $\theta: G\to R$ is the standard projection, we similarly receive that $(1-\theta)(H)=T_1\sim T_0$. But, because $T_0\oplus R_0\leq H\leq T_1\oplus R_1$, we then conclude $H\sim T_1\oplus R_1$.

Now, letting $\varphi\in\mathrm{Aut}\,G$, by assumption we have $T_1\sim\varphi(T_1)$, and thus the action of $\varphi$ on $R_1$ can be represented as the sum of some $\psi\in\mathrm{Aut}\,R$ and $\alpha\in\mathrm{Hom}(R,T)$.
Consequently, $\psi(R_1)\sim R_1$ and $\alpha(R_1)$ is finite, so that $\varphi(R_1)\sim R_1$ means $\varphi(H)\sim H$, as needed.
\fine

\begin{example}\label{0007} The next two statements are valid:

(1) Write $G=D\oplus R$, where $R$ is a torsion-free group, and $D$ is a torsion group such that each its subgroup is completely inert and, for each prime $p$, its $p$-component is a non-reduced non-zero group. Then, each subgroup in $G$ is completely inert if, and only, if $R\cong\mathbb{Z}$.

(2) If $T=\bigoplus_{p\in\Pi}T_p$ is a torsion group such that $T_p\cong\mathbb{Z}(p)$ for each $p\in\Pi$ and
$|\Pi|=\aleph_0$. Then, every direct summand of $G=\prod_{p\in\Pi}T_p$ is fully invariant, and so completely inert.
besides, in $G$ there is subgroup which is {\it not} completely inert.
\end{example}

\Pf (1) {\bf Necessity.} Suppose $H\leq R$ is a pure subgroup such that the rank of $R/H$ is exactly $1$. If
$R/H\ncong\mathbb{Z}$, then there exists a homomorphism $R\to T$ with infinite image, which is an absurd. So,
$R/H\cong\mathbb{Z}$ and hence $R$ is decomposable, which is again absurdly.

{\bf Sufficiency.} It follows at once from Proposition~\ref{0006}.

(2) Writing $G=A\oplus B$, we then arrive at $T(G)=T(A)\oplus T(B)$, and $f(T(A))=\{0\}$ for every $f\in\mathrm{Hom}(A,B)$; so, $f(A)=\{0\}$ as the quotient-group $A/T(A)$ is divisible.

Setting $x=(\dots,x_p,\dots)$, where $x_p\neq 0$, if, and only if,
$p\in\Pi$ and $X\leq G$ is such a torsion-free rank 1 subgroup that $x\in X$.

Furthermore, putting $y\in G$, $y:=(\dots,x_p',\dots)$ with $x_p'\neq 0$ for every $p\in\Pi$ such that the coset $y+T(G)$ is independent of $x+T(G)$ over $\mathbb{Q}$ in $G/T(G)$, one finds that there exists $\varphi\in\mathrm{Aut}\,G$ with $\varphi(x)=y$. That is why, a simple check shows that the subgroup $X$ is not completely inert in $G$, as stated.
\fine

An extremely difficult question is that of globally characterizing those reduced ($p$-)groups for which every their fully inert subgroup is commensurable with a fully invariant subgroup (see, e.g., \cite[Problem 2.1]{S}). The best achievements in this directions are that the class of such groups contains both the classes of totally projective groups and torsion-complete groups (see \cite{GS22}, \cite{K} and \cite[Theorem 2.2]{S}, respectively).

\medskip

That is why, it is quite logical to ask what happens in the case of characteristically inert subgroups (compare also with \cite{DK}). Concretely, one may ask to characterize those reduced ($p$-)groups whose totally (resp., completely) inert subgroups are commensurable with some fully invariant (resp., characteristic) subgroups.

\medskip

In this vein, we are now prepared to attack the following helpful observation.

\begin{proposition}\label{all} The subgroup $H$ of a group $G$ is commensurable with a characteristic subgroup of $G$ if, and only if, $H$ is uniformly completely inert in $G$.
\end{proposition}

\Pf {\bf Necessity.} It is almost evident since Proposition~\ref{3} riches us that the subgroup is uniformly completely inert if, and only if, it is uniformly characteristically inert, and each subgroup commensurable with a characteristic subgroup is uniformly characteristically inert (see \cite{DDR}).

{\bf Sufficiency.} If we assume $H$ is uniformly completely inert, then it is uniformly characteristically inert and so \cite[Corollary 1.9]{DDR} works to deduce that $H$ is commensurable with a characteristic subgroup of $G$, as required.
\fine

We, thus, now come to a major assertion which summarizes a part of the above assertions and, thereby, stimulates our further writing.

\begin{theorem}\label{fit1} Suppose $G$ is a group such that each its subgroup is commensurable with a characteristic subgroup of $G$. Then, $G$ can completely be characterized.
\end{theorem}

\Pf A consultation with Proposition~\ref{all} is a guarantor that every subgroup of $G$ is uniformly completely inert. Besides, each uniformly completely inert subgroup is manifestly completely inert. That is why, we can subsequently apply all statements starting from Lemma~\ref{0001} to Example~\ref{0007} to get the desired complete characterization after all.
\fine

Note that it cannot be derived any useful information when $G$ is a reduced group such that each its (totally inert) subgroup is commensurable with a fully invariant subgroup of $G$. Indeed, it was noted in \cite{GS} that a fully invariant subgroup $H$ of a group $G$ is totally inert in $G$ if, and only if, $f(H)$ has finite index in $H$ for all endomorphisms $f:G \to G$.

\medskip

Notice also that both Proposition~\ref{27} and Lemma~\ref{210} are common generalizations of \cite[Theorem 2.3]{S} and the same theorem appears in \cite[Theorem 3.5]{GS}.

\medskip

The following extra comments are, hopefully, worthwhile.

\begin{remark} The class of groups considered in Theorem~\ref{fit1} forms a much smaller class than the classes of groups in which every completely/characteristically inert subgroup is commensurable with a characteristic subgroup. In fact, the proper inclusion follows from the simple facts alluded to above that there are too many subgroups that are definitely {\it not} either completely nor characteristically inert. Thus, it cannot be happen that we will succeed to obtain their comprehensive descriptions at all.
\end{remark}

It was pointed out in \cite[Theorem 2.3]{S} that, if $X=\oplus_{i\in I} G_i$ is a group such that each direct summand $G_i$ is isomorphic to a fixed unbounded fully transitive $p$-group in which every fully inert subgroup is commensurable with a fully invariant subgroup, then $X$ has the same property.

\medskip

We now will expand this to the case of characteristic subgroups as follows.

\begin{theorem} Suppose that $G=\oplus_{i\in I} G_i$ is a group such that each direct summand $G_i$ is isomorphic to a fixed unbounded transitive $p$-group in which every characteristically inert subgroup is commensurable with a characteristic subgroup. Then, $G$ possesses the same property.
\end{theorem}

\Pf It entirely relies on the same arguments as in Proposition~\ref{27} and Lemma~\ref{210}, so the drop off the complete details.
\fine

Some additional things of the subject are these:

\begin{lemma}\label{ess} If $H$ is an essential subgroup in $G$ and $\phi\in\mathrm{Aut}\,G$, then $\phi(H)$ is also essential in $G$. In particular, if $H$ does not contain its proper essential subgroups, then $H$ is characteristic in $G$.
\end{lemma}

\Pf If we assume $A\cap\phi(H)=\{0\}$ for some $A\leq G$, then one verifies that $$\phi^{-1}(A\cap\phi(H))=\phi^{-1}(A)\cap H=\{0\},$$
so both $\phi^{-1}(A)$ and $A$ are zero, as required.
\fine

\begin{example} Let $G$ be a torsion-free group of finite rank, and suppose $H\leq G$ is a free essential subgroup (i.e., $r(H)=r(G)$). Then, $H$ is completely inert in $G$.
\end{example}

\Pf Letting $\phi\in\mathrm{Aut}\,G$, then Lemma~\ref{ess} employs to get that $\phi(H)$ is essential in $G$, so that $H\cap\phi(H)$ is essential in both $H$ and $\phi(H)$; in particular, $H\cap\phi(H)\cong H$ whence $H\cap\phi(H)$
has a finite index simultaneously in $H$ and $\phi(H)$. In particular, if $G$ is decomposable, then we receive additional examples of completely inert {\it not} totally inert subgroups.
\fine

We end our work with the following question which, hopefully, will stimulate a further investigation on the subject.

\medskip

\noindent{\bf Problem.} Explore those subgroups $K$ of a group $G$ such that the intersection $K\cap f(K)$ has finite index in both $K$ and $f(K)$ for all homomorphisms $f: K\to G$.

\medskip

It is apparent that these subgroups are always totally inert and strongly inert simultaneously.

\medskip
\medskip

\noindent{\bf Funding:} The scientific work of Andrey R. Chekhlov was supported by the Ministry of Science and Higher Education of Russia (agreement No. 075-02-2023-943). 

\end{document}